\documentclass[reqno]{amsart}
\usepackage{setspace}
\usepackage{graphicx}
\usepackage{amsmath}
\usepackage{amsfonts}
\usepackage{amssymb}

\begin{document}

\centerline{\huge \bf Multiple Translative Tilings in Euclidean Spaces}

\bigskip\smallskip
\centerline{\Large\bf Qi Yang and Chuanming Zong}

\vspace{0.5cm}
\centerline{\begin{minipage}{12.8cm}
{\bf Abstract.} In 1885, Fedorov discovered that a convex domain can form a lattice tiling of the Euclidean plane if and only if it is a parallelogram or a centrally symmetric hexagon. This paper proves the following results:  {\it Besides parallelograms and centrally symmetric hexagons, there is no other convex domain which can form a two-, three- or four-fold translative tiling in the Euclidean plane. However, there are convex octagons and decagons which can form five-fold translative tilings.}
\end{minipage}}

\vspace{0.4cm}

\vspace{0.8cm}
\noindent
{\Large\bf 1. Introduction}

\bigskip\noindent
In 1885, Fedorov \cite{fedo} proved that {\it a convex domain can form a lattice tiling in the plane if and only if it is a parallelogram or a centrally symmetric hexagon; a convex body can form a lattice tiling in the space if and only if it is a parallelotope, an hexagonal prism, a rhombic dodecahedron, an elongated dodecahedron, or a truncated octahedron}. As a generalized inverse problem of Fedorov's discovery, in 1900 Hilbert \cite{hilb} listed the following question as the second part of his 18th problem: {\it Whether polyhedra also exist which do not appear as fundamental regions of groups of motions, by means of which nevertheless by a suitable juxtaposition of congruent copies a complete filling up of all space is possible.} Try to verify Hilbert's problem in the plane, in 1917 Bieberbach suggested Reinhardt (see \cite{rein}) to determine all the two-dimensional convex tiles. However, to complete the list turns out to be challenging and dramatic. Over the years, the list has been successively extended by Reinhardt, Kershner, James, Rice, Stein, Mann, McLoud-Mann and Von Derau (see \cite{mann,zong14}), its completeness has been mistakenly announced several times! In 2017, M. Rao \cite{rao} announced a completeness proof based on computer checks.

Let $K$ be a convex body with (relative) interior ${\rm int}(K)$, (relative) boundary $\partial (K)$ and volume ${\rm vol}(K)$, and let $X$ be a discrete set, both in $\mathbb{E}^n$. We call $K+X$ a {\it translative tiling} of $\mathbb{E}^n$ and call $K$ a {\it translative tile} if $K+X=\mathbb{E}^n$ and the translates ${\rm int}(K)+{\bf x}_i$ are pairwise disjoint. In other words, if $K+X$ is both a packing and a covering in $\mathbb{E}^n$. In particular, we call $K+\Lambda$ a {\it lattice tiling} of $\mathbb{E}^n$ and call $K$ a {\it lattice tile} if $\Lambda $ is an $n$-dimensional lattice. Apparently, a translative tile must be a convex polytope. Usually, a lattice tile is called a {\it parallelohedron}.

To characterize all the parallelohedra in higher dimensions turns out to be very complicated. Through the works of Delone \cite{delo}, $\check{S}$togrin \cite{stog} and Engel \cite{enge}, we know that there are exact $52$ combinatorially different types of parallelohedra in $\mathbb{E}^4$. A computer classification for the five-dimensional parallelohedra was announced by Dutour Sikiri$\acute{\rm c}$, Garber, Sch$\ddot{\rm u}$rmann and Waldmann \cite{dgsw} only in 2015.

Let $\Lambda $ be an $n$-dimensional lattice. The {\it Dirichlet-Voronoi cell} of $\Lambda $ is defined by
$$C=\left\{ {\bf x}: {\bf x}\in \mathbb{E}^n,\ \| {\bf x}, {\bf o}\|\le \| {\bf x}, \Lambda \|\right\},$$
where $\| X, Y\|$ denotes the Euclidean distance between $X$ and $Y$. Clearly, $C+\Lambda $ is a lattice tiling and the Dirichlet-Voronoi cell $C$ is a parallelohedron. In 1908, Voronoi \cite{voro} made a conjecture that {\it every parallelohedron is a linear transformation image of the Dirichlet-Voronoi cell of a suitable lattice.} In $\mathbb{E}^2$, $\mathbb{E}^3$ and $\mathbb{E}^4$, this conjecture was confirmed by Delone \cite{delo} in 1929. In higher dimensions, it is still open.

To characterize the translative tiles is another fascinating problem. At the first glance, translative tilings should be more complicated than lattice tilings. However, the dramatic story had a happy end! It was shown by Minkowski \cite{mink} in 1897 that {\it every translative tile must be centrally symmetric}. In 1954, Venkov \cite{venk} proved that {\it every translative tile must be a lattice tile $($parallelohedron$)$} (see \cite{alek} for generalizations). Later, a new proof for this beautiful result was independently discovered by McMullen \cite{mcmu}.

Let $X$ be a discrete multiset in $\mathbb{E}^n$ and let $k$ be a positive integer. We call $K+X$ a {\it $k$-fold translative tiling} of $\mathbb{E}^n$ and call $K$ a {\it translative $k$-tile} if every point ${\bf x}\in \mathbb{E}^n$ belongs to at least $k$ translates of $K$ in $K+X$ and every point ${\bf x}\in \mathbb{E}^n$ belongs to at most $k$ translates of ${\rm int}(K)$ in ${\rm int}(K)+X$. In other words, $K+X$ is both a $k$-fold packing and a $k$-fold covering in $\mathbb{E}^n$. In particular, we call $K+\Lambda$ a {$k$-fold lattice tiling} of $\mathbb{E}^n$ and call $K$ a {\it lattice $k$-tile} if $\Lambda $ is an $n$-dimensional lattice. Apparently, a translative $k$-tile must be a convex polytope. In fact, similar to Minkowski's characterization, it was shown by Gravin, Robins and Shiryaev \cite{grs} that {\it a translative $k$-tile must be a centrally symmetric polytope with centrally symmetric facets.}

Multiple tilings was first investigated by Furtw\"angler \cite{furt} in 1936 as a generalization of Minkowski's conjecture on cube tilings. Let $C$ denote the $n$-dimensional unit cube. Furtw\"angler made a conjecture that {\it every $k$-fold lattice tiling $C+\Lambda$ has twin cubes. In other words, every multiple lattice tiling $C+\Lambda$ has two cubes sharing a whole facet.} In the same paper, he proved the two- and three-dimensional cases. Unfortunately, when $n\ge 4$, this beautiful conjecture was disproved by Haj\'os \cite{hajo} in 1941. In 1979, Robinson \cite{robi} determined all the integer pairs $\{ n,k\}$ for which Furtw\"angler's conjecture is false. We refer to Zong \cite{zong05,zong06} for an introduction account and a detailed account on this fascinating problem, respectively, to pages 82-84 of Gruber and Lekkerkerker \cite{grub} for some generalizations.

Let $P$ be an $n$-dimensional centrally symmetric convex polytope, let $\tau (P)$ denote the smallest integer $k$ such that $P$ can form a $k$-fold translative tiling in $\mathbb{E}^n$, and let $\tau^* (P)$ denote the smallest integer $k$ such that $P$ can form a $k$-fold lattice tiling in $\mathbb{E}^n$. For convenience, we define $\tau (P)=\infty $ if $P$ cannot form translative tiling of any multiplicity. Clearly, for every centrally symmetric convex polytope we have
$$\tau (P)\le \tau^*(P).$$

\medskip
In 1994, Bolle \cite{boll} proved that {\it every centrally symmetric lattice polygon is a lattice multiple tile}. However, little is known about the multiplicity. Let $\Lambda $ denote the two-dimensional integer lattice, and let $D'_8$ denote the octagon with vertices $(1,0)$, $(2,0)$, $(3,1)$, $(3,2)$, $(2,3)$, $(1,3)$, $(0,2)$ and $(0,1)$. As a particular example of Bolle's theorem, it was discovered by Gravin, Robins and Shiryaev \cite{grs} that {\it $D'_8+\Lambda$ is a seven-fold lattice tiling of $\mathbb{E}^2$.}

In 2017, Yang and Zong \cite{yz} studied the multiplicity of the multiple lattice tilings by proving the following results: {\it Besides parallelograms and centrally symmetric hexagons, there is no other convex domain which can form a two-, three- or four-fold lattice tiling in the Euclidean plane. However, there is a decagon $D_{10}$ which can form a five-fold lattice tilings and therefore
$$\tau^*(D_{10})=5.$$
In general, whenever $n\ge 3$, there is a non-parallelohedral polytope which can form a five-fold lattice tiling in the $n$-dimensional Euclidean space.} Afterwards, all the two-dimensional five-fold lattice tiles are characterized by Zong \cite{zong17}.

\medskip
In 2000, Kolountzakis \cite{kolo} proved that, if $D$ is a two-dimensional convex domain which is not a parallelogram and $D+X$ is a multiple tiling in $\mathbb{E}^2$, then $X$ must be a finite union of translated two-dimensional lattices. In 2013, a similar result in $\mathbb{E}^3$ was discovered by Gravin, Kolountzakis, Robins and Shiryaev \cite{gkrs}. However, up to now, no result about $\tau (P)$ is known.

\medskip
In this paper, we will investigate the multiple translative tiles by proving the following results:

\medskip\noindent
{\bf Theorem 1.} {\it Let $P_{2m}$ be a centrally symmetric convex $2m$-gon, then}
$$\tau (P_{2m})\ge \left\{\begin{array}{ll}
m-1,&\mbox{if $m$ is even,}\\
m-2,&\mbox{if $m$ is odd.}
\end{array}
\right.$$

\medskip\noindent
{\bf Theorem 2.} {\it If $D$ is a two-dimensional convex domain which is neither a parallelogram nor a centrally symmetric hexagon, then we have
$$\tau (D)\ge 5,$$
where the equality holds if $D$ is some particular centrally symmetric octagon or some particular centrally symmetric decagon.}

\vspace{0.6cm}
\noindent
{\Large\bf 2. Proof of Theorem 1}

\bigskip\noindent
Let $P_{2m}$ be a centrally symmetric convex $2m$-gon centered at the origin, with $2m$ vertices ${\bf v}_1$, ${\bf v}_2$, $\ldots$, ${\bf v}_{2m}$ enumerated in the clock order and $2m$ edges $G_1$, $G_2$, $\ldots $, $G_{2m}$, where $G_i$ is ended by ${\bf v}_i$ and ${\bf v}_{i+1}$. For convenience, we write $V=\{{\bf v}_1, {\bf v}_2, \ldots, {\bf v}_{2m}\}$ and $\Gamma=\{G_1, G_2, \ldots, G_{2m}\}$.

Assume that $P_{2m}+X$ is a $\tau (P_{2m})$-fold translative tiling in $\mathbb{E}^2$, where $X=\{{\bf x}_1, {\bf x}_2, {\bf x}_3, \ldots \}$ is a discrete multiset with ${\bf x}_1={\bf o}$. Now, let us observe the local structure of $P_{2m}+X$ at the vertices ${\bf v}\in V+X$.

Let $X^{\bf v}$ denote the subset of $X$ consisting of all points ${\bf x}_i$ such that
$${\bf v}\in \partial (P_{2m})+{\bf x}_i.$$
Since $P_{2m}+X$ is a multiple tiling, the set $X^{\bf v}$ can be divided into disjoint subsets $X^{\bf v}_1$, $X^{\bf v}_2$, $\ldots ,$ $X^{\bf v}_t$ such that the translates in $P_{2m}+X^{\bf v}_j$ can be re-enumerated as $P_{2m}+{\bf x}^j_1$, $P_{2m}+{\bf x}^j_2$, $\ldots $, $P_{2m}+{\bf x}^j_{s_j}$ satisfying the following conditions:

\medskip
\noindent
{\bf 1.} {\it ${\bf v}\in \partial (P_{2m})+{\bf x}^j_i$ holds for all $i=1, 2, \ldots, s_j.$}

\smallskip\noindent
{\bf 2.} {\it Let $\angle^j_i$ denote the inner angle of $P_{2m}+{\bf x}^j_i$ at ${\bf v}$ with two half-line edges $L^j_{i,1}$ and $L^j_{i,2}$ such that $L^j_{i,1}$, ${\bf x}^j_i-{\bf v}$ and $L^j_{i,2}$ are in clock order. Then, the inner angles join properly as
$$L^j_{i,2}=L^j_{i+1,1}$$
holds for all $i=1,$ $2,$ $\ldots ,$ $s_j$, where $L^j_{s_j+1,1}=L^j_{1,1}$.}

\medskip
For convenience, we call such a sequence $P_{2m}+{\bf x}^j_1$, $P_{2m}+{\bf x}^j_2$, $\ldots $, $P_{2m}+{\bf x}^j_{s_j}$ an {\it adjacent wheel} at ${\bf v}$. It is easy to see that
$$\sum_{i=1}^{s_j}\angle^j_i =2w_j\cdot \pi$$
hold for positive integers $w_j$. Then we define
$$\phi ({\bf v})=\sum_{j=1}^tw_j= {1\over {2\pi }}\sum_{j=1}^t\sum_{i=1}^{s_j}\angle^j_i$$
and
$$\varphi ({\bf v})=\sharp \left\{ {\bf x}_i:\ {\bf x}_i\in X,\ {\bf v}\in {\rm int}(P_{2m})+{\bf x}_i\right\}.$$

\medskip
Clearly, if $P_{2m}+X$ is a $\tau (P_{2m})$-fold translative tiling of $\mathbb{E}^2$, then
$$\tau (P_{2m})= \varphi ({\bf v})+\phi ({\bf v})\eqno (1)$$
holds for all ${\bf v}\in V+X$.

\medskip
Now, we introduce a basic lemma which is useful in the studying of $\tau (P_{2m}).$

\medskip\noindent
{\bf Lemma 1.} {\it Assume that $P_{2m}$ is a centrally symmetric convex $2m$-gon centered at the origin and $P_{2m}+X$ is a $\tau (P_{2m})$-fold translative tiling of the plane, where $m\ge 4$. If ${\bf v}\in V+X$ is a vertex and $G\in \Gamma +X$ is an edge with ${\bf v}$ as one of its two ends, then there are at least $\lceil (m-3)/2\rceil $ different translates $P_{2m}+{\bf x}_i$ satisfying both
$${\bf v}\in \partial (P_{2m})+{\bf x}_i$$
and}
$$G\setminus \{ {\bf v}\}\subset {\rm int}(P_{2m})+{\bf x}_i.$$

\smallskip
\noindent
{\bf Proof.} Since adjacent wheels are circular, without loss of generality, let $P_{2m}+{\bf x}_1$, $P_{2m}+{\bf x}_2$, $\ldots$, $P_{2m}+{\bf x}_s$ be an adjacent wheel at ${\bf v}$ such that $G$ is the first edge appearing in the wheel and let $\angle_i$ denote the inner angle of $P_{2m}+{\bf x}_i$ at the vertex ${\bf v}$.

Let $n$ denote the smallest index such that
$$\sum_{i=1}^n\angle_i=\omega \cdot \pi \eqno(2)$$
holds with some positive integer $\omega $. If $\angle_j$ and $\angle_{j+k}$ are two opposite angles of $P_{2m}$ appearing in the angle sequence with $1\le j<j+k\le n$, it is easy to see that
$$\sum_{i=0}^{k-1}\angle_{j+i}=\omega'\cdot \pi$$
holds with a positive integer $\omega'$ and $\omega\ge \omega'.$ Therefore, to estimate $\omega $, we may assume that the angle sequence $\angle_1$, $\angle_2$, $\ldots $, $\angle_n$ has no opposite angle pair of $P_{2m}$.

\begin{figure}[!ht]
\centering
\includegraphics[scale=0.55]{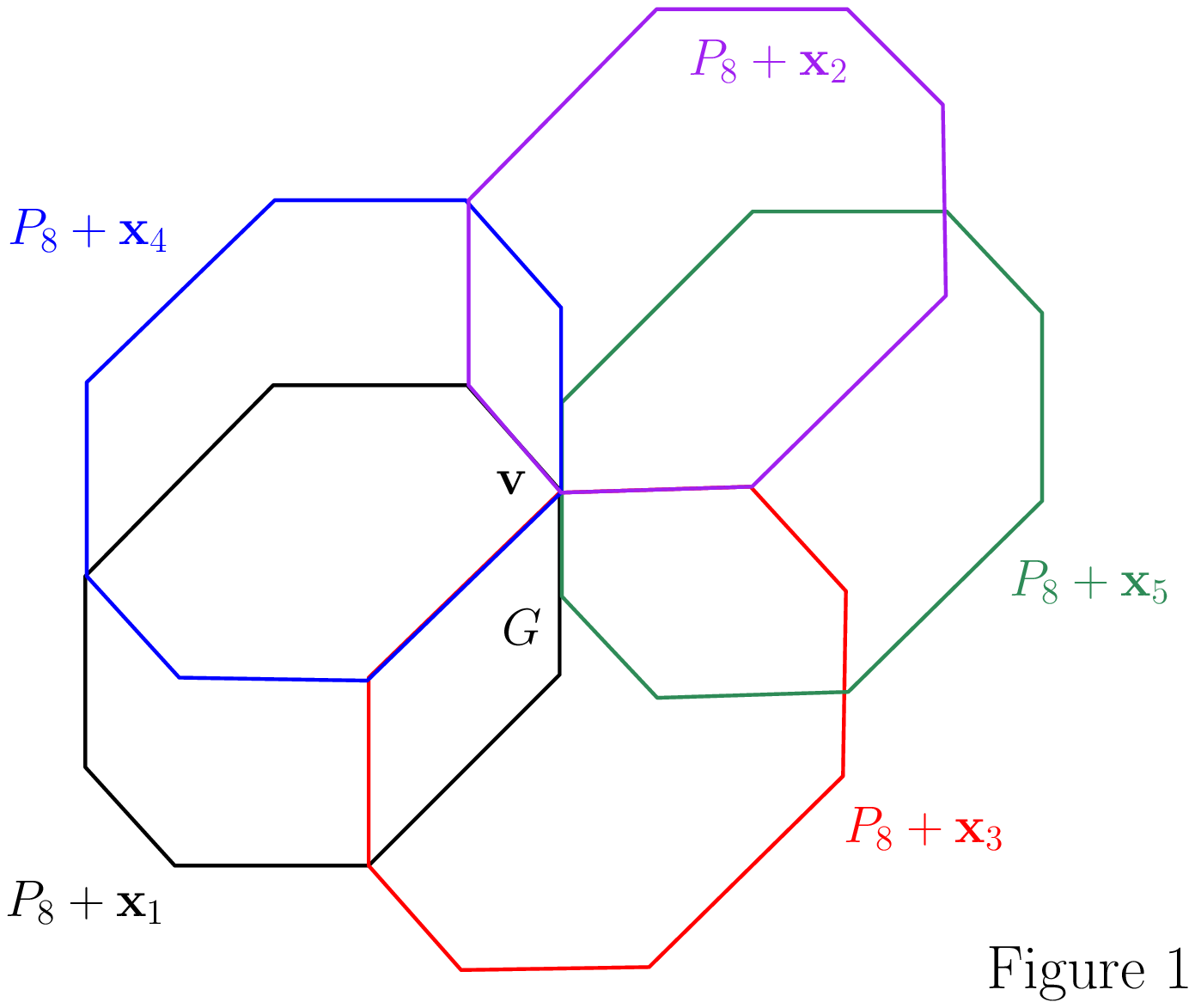}
\end{figure}

\medskip
Clearly, $\angle_i=\pi$ if and only if ${\bf v}$ is a relative interior point of an edge of $P_{2m}+{\bf x}_i$ (such as $\angle_5$ in Figure 1) and therefore
$$\sum_{i=1}^n\angle_i<n\cdot \pi.\eqno(3)$$
On the other hand, if $\ell$ of the $n$ angles are $\pi $ and $n-\ell<m$, then we have
$$\sum_{i=1}^n\angle_i>\ell\cdot \pi +(m-1)\cdot\pi -(m-n+\ell )\cdot\pi =(n-1)\cdot \pi,\eqno (4)$$
which together with (3) contradicts (2). Therefore, to avoid the contradiction, we must have
$$n-\ell =m$$
and each pair of the opposite angles of $P_{2m}$ has a representative in the sequence $\angle_1$, $\angle_2$, $\ldots ,$ $\angle_n$. Consequently, we have
$$\sum_{i=1}^n\angle_i \ge {{(2m-2)\cdot \pi }\over 2}=(m-1)\cdot \pi.\eqno (5)$$
Therefore $G\setminus \{{\bf v}\}$ is covered by at least
$$\left\lceil {{m-1}\over 2}\right\rceil -1=\left\lceil {{m-3}\over 2}\right\rceil $$
of the $s$ translates ${\rm int}(P_{2m})+{\bf x}_i$. Lemma 1 is proved.\hfill{$\Box$}

\bigskip\noindent
{\bf Proof of Theorem 1.} Assume that $P_{2m}+X$ is a $\tau (P_{2m})$-fold translative tiling in the Euclidean plane and assume that ${\bf v}\in V+X$. Then it follows by Lemma 1 that
$$\varphi({\bf v})\ge \left\lceil {{m-3}\over 2}\right\rceil .\eqno(6)$$
Let $P_{2m}+{\bf x}_1$, $P_{2m}+{\bf x}_2$, $\ldots $, $P_{2m}+{\bf x}_s$ be an adjacent wheel at ${\bf v}$ and let
$\angle_1$, $\angle_2,$ $\ldots ,$ $\angle_s$ be the corresponding angle sequence. By (5) we have
$$\phi ({\bf v})\ge {1\over {2\pi}}\sum_{i=1}^s\angle_i\ge \left\lceil {{m-1}\over 2}\right\rceil.\eqno(7)$$

Then, it follows by (1), (6) and (7) that
$$\tau (P_{2m})\ge \left\lceil {{m-3}\over 2}\right\rceil +\left\lceil {{m-1}\over 2}\right\rceil =\left\{
\begin{array}{ll}
m-1,&\mbox{if $m$ is even,}\\
m-2,&\mbox{if $m$ is odd.}
\end{array}
\right.$$
Theorem 1 is proved. \hfill{$\Box$}

\vspace{0.8cm}
\noindent
{\Large\bf 3. Proof of Theorem 2}

\bigskip\noindent
First, let's introduce a formula for $\phi ({\bf v})$ which will be useful.

\medskip\noindent
{\bf Lemma 2.} {\it Assume that $P_{2m}$ is a centrally symmetric convex $2m$-gon centered at the origin, $P_{2m}+X$ is a translative multiple tiling of the plane, and ${\bf v}\in V+X$. Then we have
$$\phi ({\bf v})=\kappa\cdot {{m-1}\over 2}+\ell \cdot {1\over 2},$$
where $\kappa $ is a positive integer and $\ell $ is a nonnegative integer. In fact, $\ell $ is the number of the edges in $\Gamma +X$ which take ${\bf v}$ as an interior point.}

\medskip
\noindent
{\bf Proof.} Assume that $P_{2m}+{\bf x}_1$, $P_{2m}+{\bf x}_2$, $\ldots,$ $P_{2m}+{\bf x}_s$ is an adjacent wheel at ${\bf v}$ and let $\angle _i$ denote the inner angle of $P_{2m}+{\bf x}_i$ at ${\bf v}$. Of course, we have $\angle_i=\pi $ if ${\bf v}$ is not a vertex of $P_{2m}+{\bf x}_i$.

Assume that $\angle_1<\pi $ and let $n$ to be the smallest index such that
$$\sum_{i=1}^n\angle_i=\omega \pi\eqno(8)$$
holds with a positive integer $\omega $. We proceed to show that each pair of the opposite angles of $P_{2m}$ has one and only one representative in
$\angle_1$, $\angle_2$, $\ldots ,$ $\angle_n$.

If, on the contrary, $\angle_j$ and $\angle_{j+k}$ are two of these $n$ angles, $\angle_j<\pi $, which are either identical or opposite. Then, it is easy to see that
$$\sum_{i=0}^{k-1}\angle_{j+i} =\omega' \pi \eqno(9)$$
holds with a positive integer $\omega'$. For convenience, we assume that $\angle_j$, $\angle_{j+1}$, $\ldots ,$ $\angle_{j+k-1}$ have neither identical nor opposite pair. Then, by repeating the argument between (2) and (5) in the proof of Lemma 1, one can deduce that each pair of the opposite angles of $P_{2m}$ has one and only one representative in $\angle_j$, $\angle_{j+1}$, $\ldots ,$ $\angle_{j+k-1}$. Consequently, one of these $k$ angles is either identical or opposite to $\angle_1$, which contradicts the minimum assumption on $n$ and $\omega$.

Then, applying the argument between (2) and (5) to $\angle_1$, $\angle_2$, $\ldots ,$ $\angle_n$, it can be deduced that
$$\sum_{i=1}^n\angle_i=(m-1) \pi +\ell_1 \pi, \eqno(10)$$
where $\ell_1$ is the number of the $\pi $ angles in $\angle_1$, $\angle_2$, $\ldots ,$ $\angle_n$. In fact, it is $n-m$.

By repeating this process to $\angle_{n+1}$, $\angle_{n+2}$, $\ldots,$ $\angle_s$ if necessary, it follows that
$$\sum_{i=1}^s\angle_i=\kappa' (m-1) \pi +\ell' \pi \eqno(11)$$
and therefore
$$\phi({\bf v})={1\over {2\pi }}\sum \sum_{i=1}^s\angle_i=\kappa\cdot {{m-1}\over 2}+\ell \cdot {1\over 2},\eqno (12)$$
where the first sum is over all adjacent wheels at ${\bf v}$, $\kappa'$ and $\kappa $ are suitable positive integers, $\ell'$ and $\ell $ are suitable nonnegative integers. In fact, $\ell$ is the number of the edges which take ${\bf v}$ as an interior point.

Lemma 2 is proved.\hfill{$\Box$}

\medskip\noindent
{\bf Lemma 3.} {\it Let $P_8$ be a centrally symmetric convex octagon, then we have
$$\tau (P_8)\ge 5,$$
where the equality holds at some particular octagon.}

\medskip
\noindent
{\bf Proof.} First of all, it follows from Lemma 1 that
$$\varphi ({\bf v})\ge \left\lceil {{4-3}\over 2}\right\rceil =1\eqno(13)$$
holds for all ${\bf v}\in V+X.$  On the other hand, by Lemma 2 we have
$$\phi ({\bf v})=\kappa \cdot {3\over 2}+\ell \cdot {1\over 2}\ge 2,\eqno(14)$$
 where $\kappa $ is a positive integer and $\ell $ is a nonnegative integer. Thus, to prove the lemma it is sufficient to deal with the following three cases:

\medskip\noindent
{\bf Case 1.} {\it $\phi ({\bf v})\ge 4$ holds for a vertex ${\bf v}\in V+X$}. Then, by (1) and (13) we get
$$\tau (P_8)=\varphi ({\bf v})+\phi ({\bf v})\ge 5.\eqno(15)$$

\medskip\noindent
{\bf Case 2.} {\it $\phi ({\bf v})=3$ holds for a vertex ${\bf v}\in V+X$}. If ${\bf v}\in {\rm int}(G)$
holds for some $G\in \Gamma +X$, then it follows by Lemma 1 that
$$\varphi ({\bf v})\ge 2$$
and therefore
$$\tau (P_8)=\varphi ({\bf v})+\phi ({\bf v})\ge 5.\eqno(16)$$

\begin{figure}[!ht]
\centering
\includegraphics[scale=0.55]{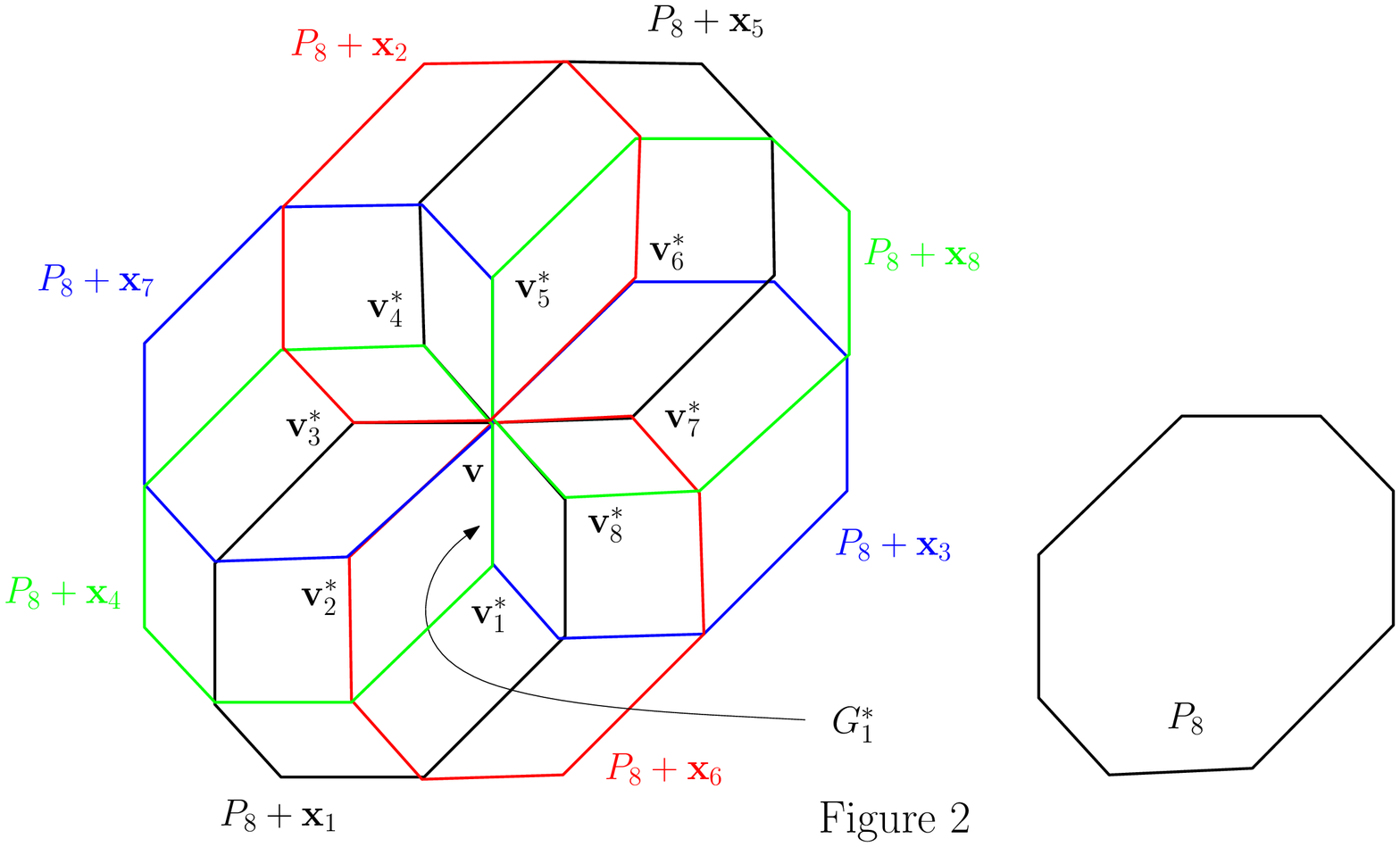}
\end{figure}

If ${\bf v}$ is a vertex for all $P_8+{\bf x}$, ${\bf x}\in X^{\bf v}$, then one can deduce that $P_8+X^{\bf v}$ is an adjacent wheel of eight translates $P_8+{\bf x}_1$, $P_8+{\bf x}_2$, $\ldots ,$ $P_8+{\bf x}_8$ as shown by Figure 2. Let ${\bf v}^*_i$ be the vertex connecting to ${\bf v}$ by edge $G^*_i$ as shown by Figure 2. By Lemma 1, for each ${\bf v}^*_i$ there is a ${\bf y}_i\in X$ satisfying both
$${\bf v}\in {\rm int}(P_8)+{\bf y}_i$$
and
$${\bf v}^*_i\in \partial (P_8)+{\bf y}_i.$$
If ${\bf y}_1={\bf y}_2=\ldots ={\bf y}_8$, since $(P_8+{\bf x})\cap (P_8+{\bf y})$ is always centrally symmetric, one can deduce that both
$(P_8+{\bf y}_1)\cap (P_8+{\bf x}_3)$ and $(P_8+{\bf y}_1)\cap (P_8+{\bf x}_4)$ are parallelograms,
$$G^*_1+{\bf v}^*_6-{\bf v}\in \Gamma +X$$
and
$$G^*_1+{\bf v}^*_4-{\bf v}\in \Gamma +X.$$
By symmetry, it follows that $P_8$ is an hexagon, which contradicts the assumption that it is an octagon. Therefore, we have
$$\varphi ({\bf v})\ge 2$$
and
$$\tau (P_8)=\varphi ({\bf v})+\phi ({\bf v})\ge 5.\eqno(17)$$

\medskip\noindent
{\bf Case 3.} {\it $\phi ({\bf v})=2$ holds for all vertices ${\bf v}\in V+X$}. First, by (14) it follows that $\phi ({\bf v})=2$ if and only if $P_8+X^{\bf v}$ is an adjacent wheel of five translates. By re-enumeration we may assume that $\angle_1$, $\angle_2$, $\angle_3$ and $\angle_4$ are inner angles of $P_8$ and $\angle_5=\pi$, as shown by Figure 3. For intuitive reason, guaranteed by linear transformation, we assume that the edges $G_1$ and $G_3$ of $P_8$ are horizontal and vertical, respectively.

\begin{figure}[!ht]
\centering
\includegraphics[scale=0.55]{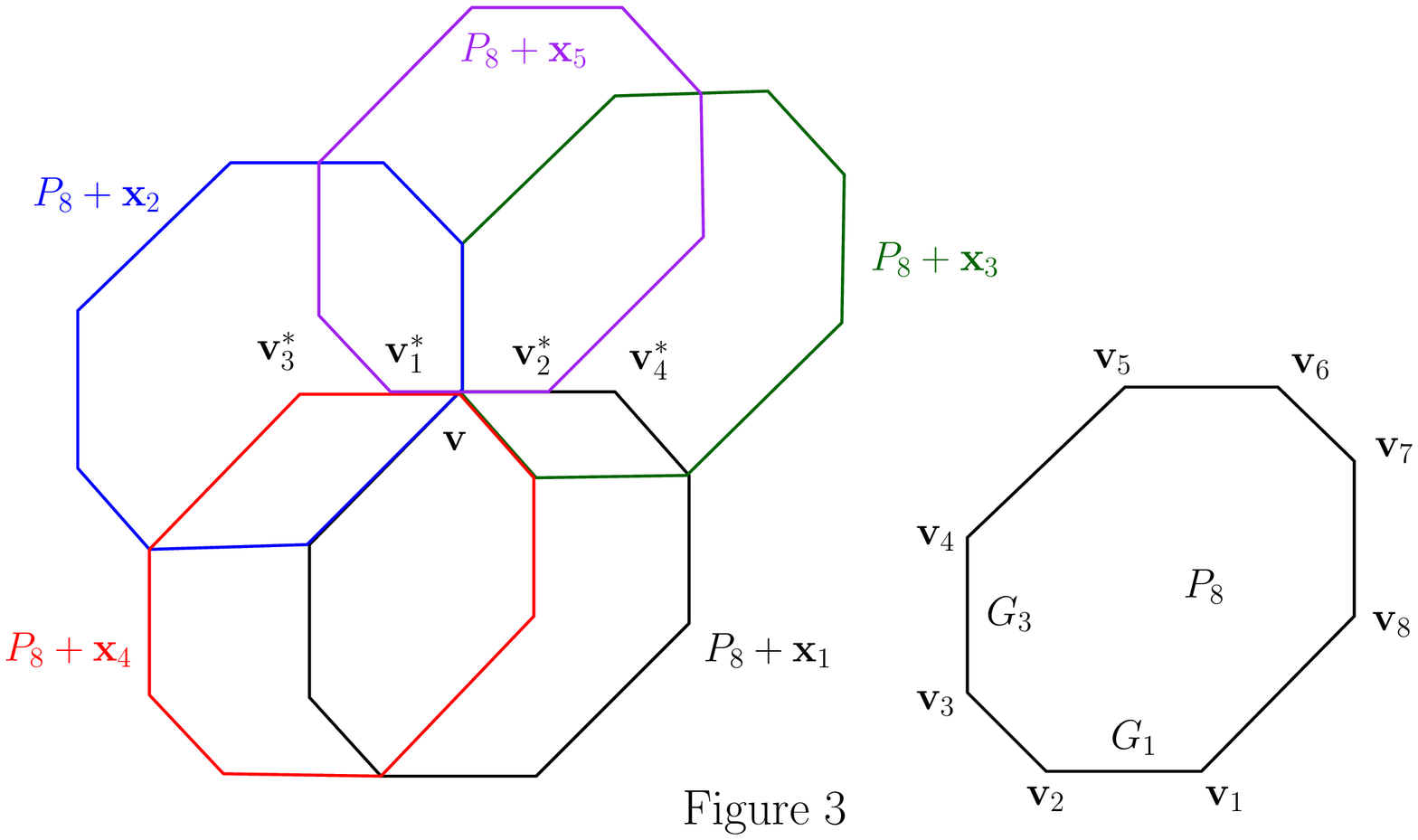}
\end{figure}

Let $G^*_1$ denote the edge of $P_8+{\bf x}_5$ such that ${\bf v}\in {\rm int}(G^*_1)$ with two ends ${\bf v}^*_1$ and ${\bf v}^*_2$, let $L$ denote the straight line determined by ${\bf v}^*_1$ and ${\bf v}^*_2$, let $G^*_3$ denote the edge of $P_8+{\bf x}_4$ lying on $L$ with ends ${\bf v}$ and ${\bf v}^*_3$, and let $G^*_4$ denote the edge of $P_8+{\bf x}_1$ lying on $L$ with ends ${\bf v}$ and ${\bf v}^*_4$.

\begin{figure}[!ht]
\centering
\includegraphics[scale=0.55]{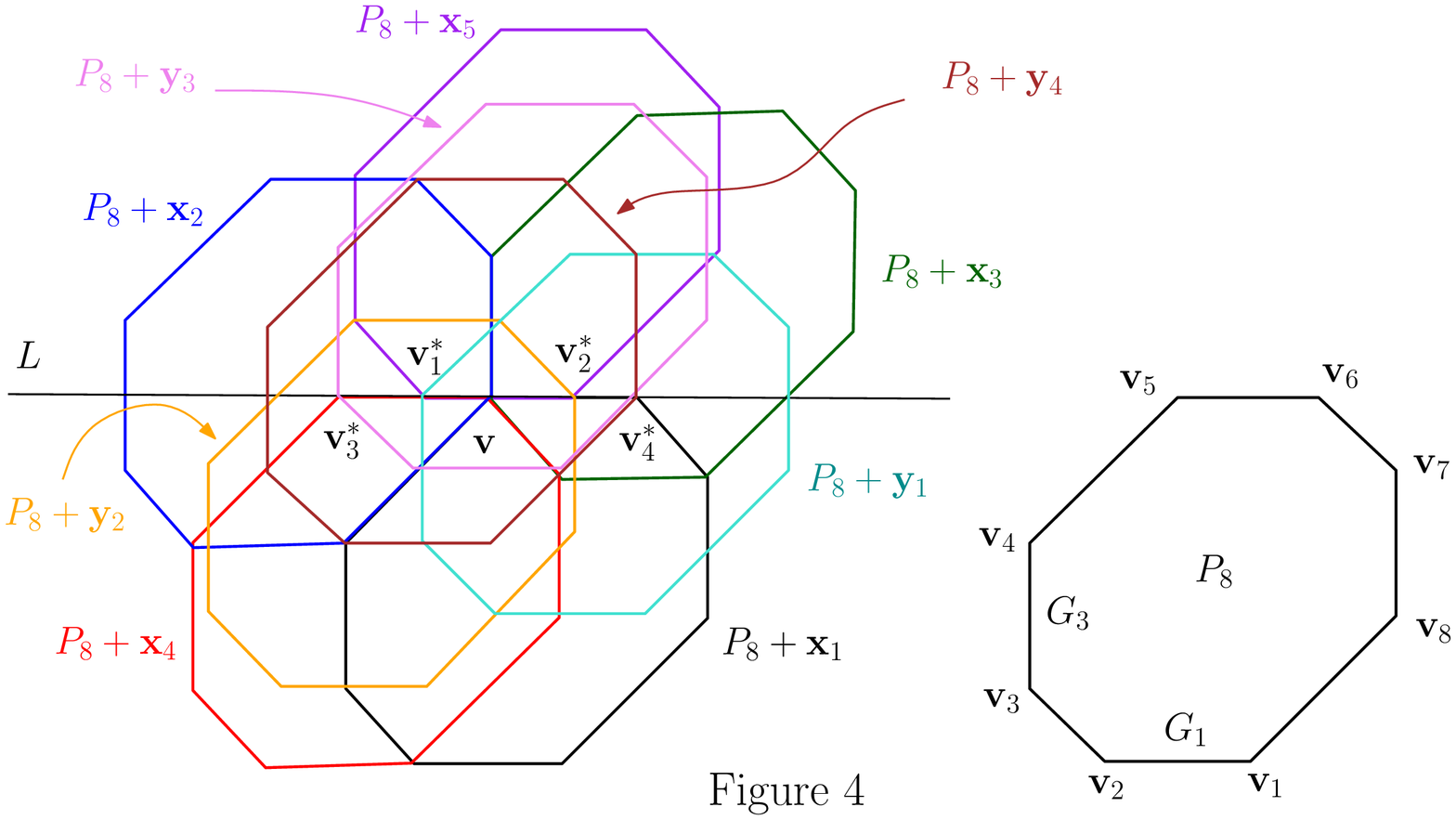}
\end{figure}

Since Figure 3 is the only configuration for $\phi ({\bf v})=2$ and $$\phi ({\bf v}^*_1)=\phi ({\bf v}^*_2)=\phi ({\bf v}^*_3)=\phi ({\bf v}^*_4)=2,\eqno(18)$$
as shown by Figure 4 the set $X$ has four points ${\bf y}_1$, ${\bf y}_2$, ${\bf y}_3$ and ${\bf y}_4$ such that
$${\bf v}^*_1={\bf v}_4+{\bf y}_1,\quad {\bf v}\in {\rm int}(P_8)+{\bf y}_1,\eqno(19)$$
$${\bf v}^*_2={\bf v}_7+{\bf y}_2,\quad {\bf v}\in {\rm int}(P_8)+{\bf y}_2,\eqno(20)$$
$${\bf v}^*_3={\bf v}_3+{\bf y}_3,\quad {\bf v}\in {\rm int}(P_8)+{\bf y}_3\eqno(21)$$
and
$${\bf v}^*_4={\bf v}_8+{\bf y}_4,\quad {\bf v}\in {\rm int}(P_8)+{\bf y}_4.\eqno(22)$$

Clearly, by the convexity of $P_8$ we have ${\bf y}_1\not= {\bf y}_2$, ${\bf y}_1\not= {\bf y}_3$ and ${\bf y}_2\not= {\bf y}_4$. For convenience, we write ${\bf v}_i=(x_i, y_i)$. If ${\bf y}_2={\bf y}_3$, then it follows from (20) and (21) that
$$y_3=y_7.\eqno(23)$$
If ${\bf y}_1={\bf y}_4$, then it follows from (19) and (22) that
$$y_4=y_8.\eqno(24)$$
It is obvious that (23) and (24) cannot hold simultaneously. Therefore, we get
$$\varphi ({\bf v})\ge 3$$
and
$$\tau (P_8)=\phi ({\bf v})+\varphi ({\bf v})\ge 5.\eqno(25)$$

\medskip
As a conclusion of these three cases, we have proved that
$$\tau (P_8)\ge 5\eqno (26)$$
holds for every centrally symmetric octagon.

\begin{figure}[!ht]
\centering
\includegraphics[scale=0.8]{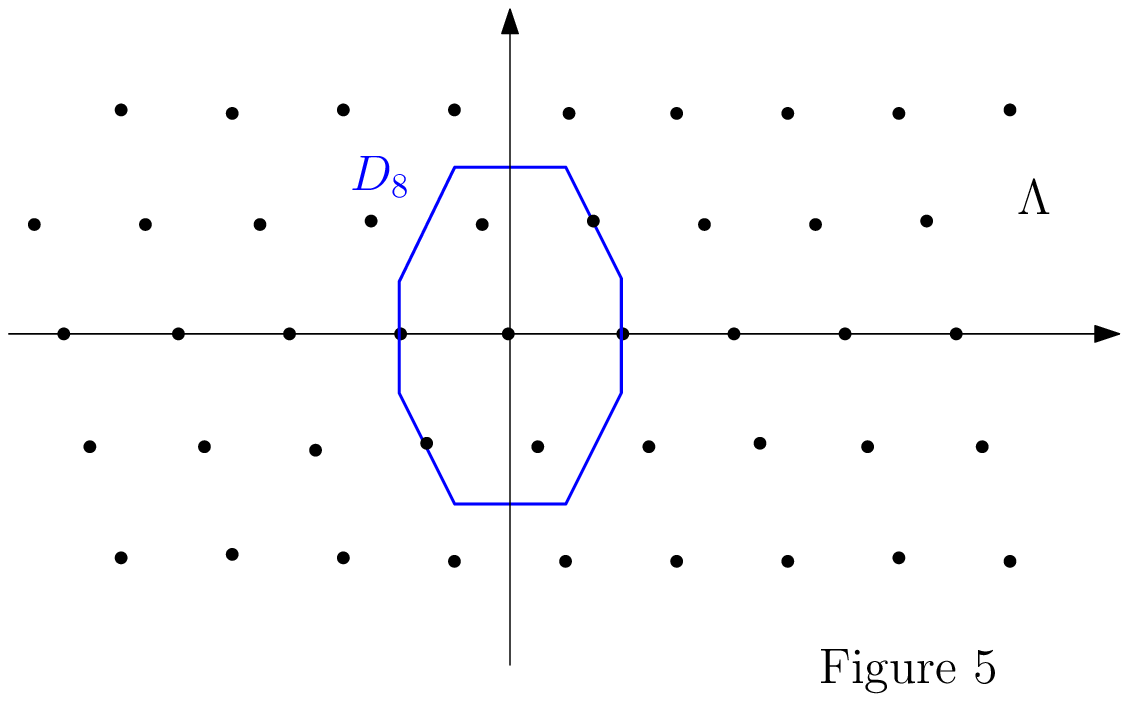}
\end{figure}

\medskip
Let $D_8$ denote the octagon with vertices ${\bf v}_1=(-1,3)$, ${\bf v}_2=(1,3)$, ${\bf v}_3=(2,1)$, ${\bf v}_4=(2,-1)$,
${\bf v}_5=-{\bf v}_1$, ${\bf v}_6=-{\bf v}_2$, ${\bf v}_7=-{\bf v}_3$ and ${\bf v}_8=-{\bf v}_4$, and let $\Lambda$ denote the lattice generated by
${\bf u}_1=(2,0)$ and ${\bf u}_2=({3\over 2},2)$ (as shown by Figure 5). Namely, 
$$\Lambda =\left\{ z_1{\bf u}_1+z_2{\bf u}_2:\ z_i\in \mathbb{Z}\right\}.\eqno(27)$$
It can be verified that $D_8+\Lambda$ is a five-fold translative tiling. Therefore, we have
$$\tau (D_8)=5.\eqno(28)$$

Lemma 3 is proved.\hfill{$\Box$}

\bigskip\noindent
{\bf Lemma 4.} {\it Let $P_{10}$ be a centrally symmetric convex decagon, then we have
$$\tau (P_{10})\ge 5,$$
where the equality holds at some particular decagon.}

\medskip
\noindent
{\bf Proof.} First of all, it follows from Lemma 1 that
$$\varphi ({\bf v})\ge \left\lceil {{5-3}\over 2}\right\rceil =1\eqno (29)$$
holds for all ${\bf v}\in V+X.$ On the other hand, by Lemma 2 we have
$$\phi ({\bf v})=\kappa \cdot 2+\ell \cdot {1\over 2},\eqno(30)$$
 where $\kappa$ is a positive integer and $\ell$ is the number of the edges which take ${\bf v}$ as an interior point. Thus, to prove the lemma it is sufficient to deal with the following three cases:

\medskip\noindent
{\bf Case 1.} {\it $\phi ({\bf v})\ge 4$ holds for a vertex ${\bf v}\in V+X$}.
Then, by (1) and (29) we have
$$\tau (P_{10})=\varphi ({\bf v})+\phi ({\bf v})\ge 5.\eqno(31)$$

\medskip\noindent
{\bf Case 2.} {\it $\phi ({\bf v})=3$ holds for a vertex ${\bf v}\in V+X$}. It follows by (30) that $\ell\not= 0$ and thus
$${\bf v}\in {\rm int}(G)$$
holds for some $G\in \Gamma +X$. Let ${\bf v}'$ and ${\bf v}^*$ be the two ends of $G$. It follows by Lemma 1 that $P_{10}+X$ has two translates $P_{10}+{\bf x}'$ and $P_{10}+{\bf x}^*$ satisfying
$${\bf v}\in {\rm int}(P_{10})+{\bf x}',\quad {\bf v}'\in \partial (P_{10})+{\bf x}'$$
and
$${\bf v}\in {\rm int}(P_{10})+{\bf x}^*,\quad {\bf v}^*\in \partial (P_{10})+{\bf x}^*.$$
By the convexity of $P_{10}$ it is easy to see that ${\bf x}'\not= {\bf x}^*$ and therefore
$$\varphi ({\bf v})\ge 2.$$
Then, it follows by (1) that
$$\tau (P_{10})=\varphi ({\bf v})+\phi ({\bf v})\ge 5.\eqno(32)$$

\smallskip\noindent
{\bf Case 3.} {\it $\phi ({\bf v})=2$ holds for a vertex ${\bf v}\in V+X$}. Then five translates $P_{10}+{\bf x}_1$, $P_{10}+{\bf x}_2$, $\ldots $, $P_{10}+{\bf x}_5$ are successively adjacent as shown by Figure 5. Let ${\bf v}^*_1$, ${\bf v}^*_2$, $\ldots,$ ${\bf v}^*_5$ be the five neighbor vertices of ${\bf v}$ as shown by Figure 6. By Lemma 1, for each of these five vertices ${\bf v}^*_i$ there is a ${\bf y}_i\in X$ satisfying both
$${\bf v}\in {\rm int}(P_{10})+{\bf y}_i$$
and
$${\bf v}^*_i\in \partial (P_{10})+{\bf y}_i.$$

\begin{figure}[!ht]
\centering
\includegraphics[scale=0.5]{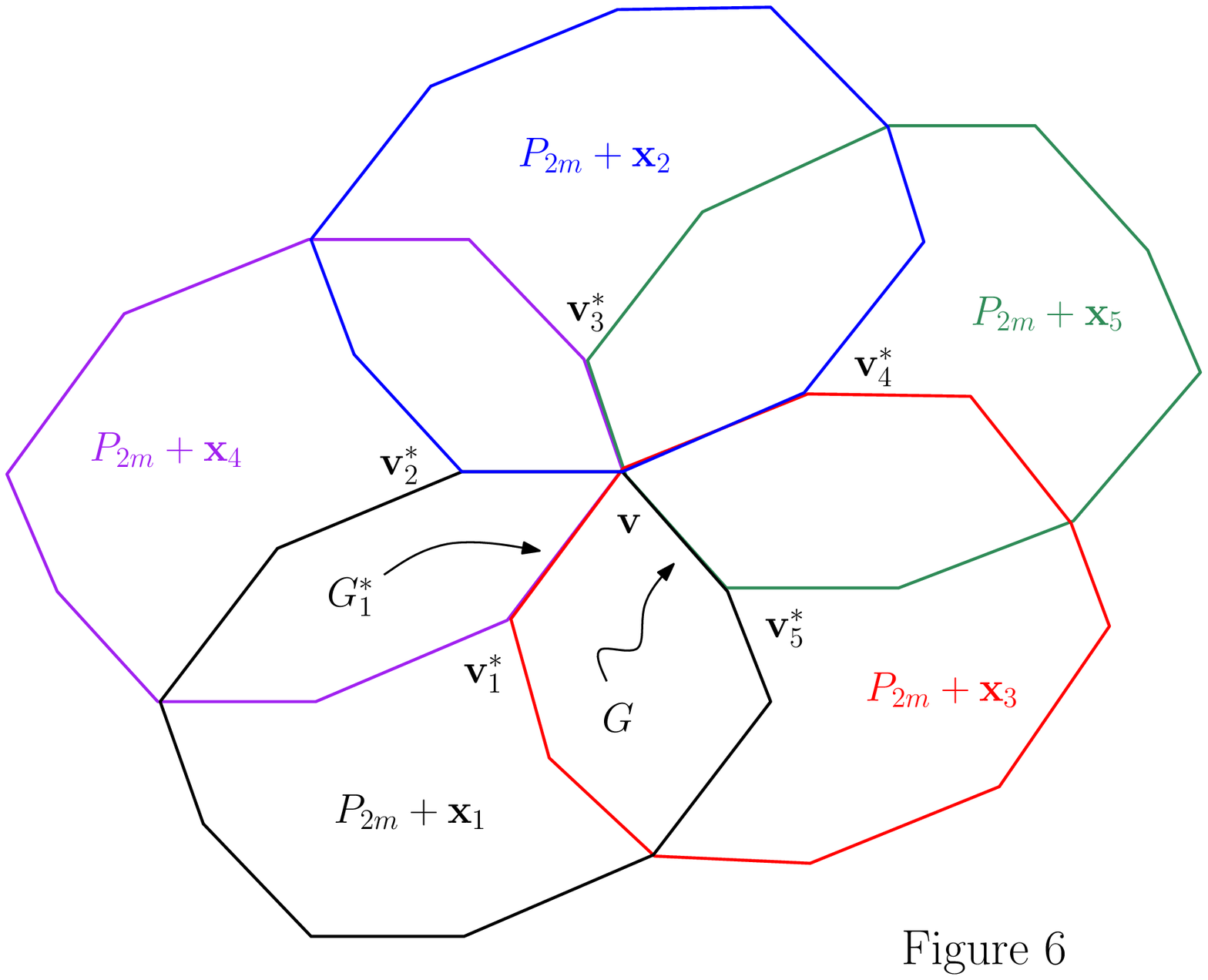}
\end{figure}

It is known that $(D+{\bf x})\cap (D+{\bf y})$ is centrally symmetric for all ${\bf x}$ and ${\bf y}$ whenever $D$ is centrally symmetric. If three of the five neighbor vertices, say ${\bf v}^*_1$, ${\bf v}^*_2$ and ${\bf v}^*_3$, are successively on the boundary of $P_{10}+{\bf y}_1$, then $(P_{10}+{\bf y}_1)\cap (P_{10}+{\bf x}_4)$ is a parallelogram. Consequently, letting $G^*_1$ denote the edge with ends ${\bf v}$ and ${\bf v}^*_1$, we have
$$G^*_1+{\bf v}^*_3-{\bf v}\in \Gamma +X.$$
This means that two opposite edges joining at ${\bf v}^*_3$ and therefore
$$\phi ({\bf v}^*_3)\ge 3.$$
Then, by the two previous cases, we have
$$\tau (P_{10})=\varphi ({\bf v}^*_3)+ \phi ({\bf v}^*_3)\ge 5.\eqno(33)$$

If three of the five neighbor vertices, say ${\bf v}^*_1$, ${\bf v}^*_3$ and ${\bf v}^*_4$, are on the boundary of $P_{10}+{\bf y}_1$ but not successively, then both $(P_{10}+{\bf y}_1)\cap (P_{10}+{\bf x}_4)$ and $(P_{10}+{\bf y}_1)\cap (P_{10}+{\bf x}_3)$ are parallelograms. Consequently, we have both
$$G^*_1+{\bf v}^*_3-{\bf v}\in \Gamma +X$$
and
$$G^*_1+{\bf v}^*_4-{\bf v}\in \Gamma +X.$$
Then, by symmetry one can deduce that $P_{10}$ is an hexagon, which contradicts the assumption that $P_{10}$ is a decagon.

As a conclusion, for any of these ${\bf y}_i$ there are at most two of the five vertices ${\bf v}^*_1$, ${\bf v}^*_2$, $\ldots ,$ ${\bf v}^*_5$
can appear on the boundary of $P_{10}+{\bf y}_i.$ Therefore, we have
$$\varphi ({\bf v})\ge 3$$
and
$$\tau (P_{10})=\varphi ({\bf v})+ \phi ({\bf v}) \ge 5.\eqno(34)$$

As a conclusion of these three cases, we have proved that
$$\tau (P_{10})\ge 5\eqno(35)$$
holds for every centrally symmetric decagon.

On the other hand, let $\Lambda$ denote the integer lattice $\mathbb{Z}^2$ and let $D_{10}$ denote the decagon with vertices ${\bf v}_1=(-{3\over 5}, -{5\over 4}),$ ${\bf v}_2=({3\over 5}, -{3\over 4}),$ ${\bf v}_3=({7\over 5}, -{1\over 4}),$ ${\bf v}_4=({8\over 5}, {1\over 4}),$
${\bf v}_5=({7\over 5}, {3\over 4}),$ ${\bf v}_6=-{\bf v}_1$, ${\bf v}_7=-{\bf v}_2$, ${\bf v}_8=-{\bf v}_3$, ${\bf v}_9=-{\bf v}_4$ and ${\bf v}_{10}=-{\bf v}_5$, it was discovered by Yang and Zong \cite{yz} that $D_{10}+\Lambda$ is a five-fold translative tiling. Thus we have
$$\tau (D_{10})=5.\eqno(36)$$

Lemma 4 is proved. \hfill{$\Box$}

\medskip
\noindent
{\bf Lemma 5.} {\it Let $P_{12}$ be a centrally symmetric convex dodecagon, then we have}
$$\tau (P_{12})\ge 6.$$

\medskip
\noindent
{\bf Proof.} First of all, it follows by Lemma 1 that
$$\varphi ({\bf v})\ge \left\lceil {{6-3}\over 2}\right\rceil =2 \eqno(37)$$
holds for all ${\bf v}\in V+X$. On the other hand, by studying the angle sums of the adjacent wheels at ${\bf v}$ one can deduce that
$$\phi ({\bf v})\ge  \left\lceil {{6-1}\over 2}\right\rceil =3.\eqno(38)$$
Thus, to show the lemma it is sufficient to deal with the following two cases.

\medskip\noindent
{\bf Case 1.} {\it $\phi({\bf v})\ge 4$ holds for a vertex ${\bf v}\in V+X$.} Then it follows by (1) and (37) that
$$\tau (P_{12})=\varphi ({\bf v})+\phi ({\bf v})\ge 6.\eqno(39)$$

\medskip\noindent
{\bf Case 2.} {\it $\phi({\bf v})=3$ holds for a vertex ${\bf v}\in V+X$.} Assume that $P_{12}+{\bf x}_1$, $P_{12}+{\bf x}_2$, $\ldots $, $P_{12}+{\bf x}_s$ is an adjacent wheel at ${\bf v}$. By studying the corresponding angle sum, it can be deduced that there is a $G\in \Gamma +X$ such that
$${\bf v}\in {\rm int}(G).$$

Let ${\bf v}'$ and ${\bf v}^*$ denote the two ends of $G$. By Lemma 1 and the convexity of $P_{12}$ it follows that $X$ has four different points ${\bf y}'_1$, ${\bf y}'_2$,
${\bf y}^*_1$ and ${\bf y}^*_2$ satisfying
$${\bf v}'\in \partial (P_{12})+{\bf y}'_i,\quad i=1,\ 2,$$
$${\bf v}^*\in \partial (P_{12})+{\bf y}^*_i,\quad i=1,\ 2,$$
$${\bf v}\in {\rm int}(P_{12})+{\bf y}'_i,\quad i=1,\ 2$$
and
$${\bf v}\in {\rm int}(P_{12})+{\bf y}^*_i,\quad i=1,\ 2.$$
Consequently, we have
$$\varphi ({\bf v})\ge 4$$
and therefore
$$\tau (P_{12})=\varphi ({\bf v})+\phi ({\bf v})\ge 7.\eqno(40)$$

As a conclusion of these two cases that
$$\tau (P_{12})\ge 6\eqno(41)$$
holds for every centrally symmetric dodecagon. Lemma 5 is proved. \hfill{$\Box$}

\medskip
\noindent
{\bf Lemma 6 (Zong \cite{zong17}).} {\it Let $P_{14}$ be a centrally symmetric convex tetradecagon, then we have}
$$\tau (P_{14})\ge 6.$$

\medskip\noindent
{\bf Proof of Theorem 2.} The theorem follows from Lemmas 3, 4, 5, 6 and Theorem 1. \hfill{$\Box$}

\vspace{0.8cm}\noindent
{\bf Acknowledgements.} For helpful comments and suggestions, the authors are grateful to Professor S. Robins and Professor G. M. Ziegler. This work is supported by 973 Program 2013CB834201 and the Chang Jiang Scholars Program of China.

\bibliographystyle{amsplain}

\begin{thebibliography}{99}
\bibitem{alek}A. D. Aleksandrov, On tiling space by polytopes, {\it Vestnik Leningrad Univ. Ser. Mat. Fiz. Him}. {\bf 9} (1954), 33-43.
\bibitem{boll}U. Bolle, On multiple tiles in $R^2$, {\it Intuitive Geometry,} Colloq. Math. Soc. J. Bolyai {\bf 63}, North-Holland, Amsterdam, 1994.
\bibitem{delo}B. N. Delone, Sur la partition reguli$\grave{e}$re de l'espace $\grave{a}$ $4$ dimensions I, II,
{\it Izv. Akad. Nauk SSSR, Ser. VII} (1929), 79-110; 147-164.
\bibitem{dgsw}M. Dutour Sikiri$\acute{\rm c}$, A. Garber, A. Sch$\ddot{\rm u}$rmann and C. Waldmann, The complete classification of five-dimensional Dirichlet-Voronoi polyhedra of translational lattices,  {\it Acta Crystallographica} {\bf A72} (2016), 673-683.
\bibitem{enge}P. Engel, On the symmetry classification of the four-dimensional parallelohedra, {\it Z. Kristallographie} {\bf 200} (1992), 199-213.
\bibitem{fedo}E. S. Fedorov, Elements of the study of figures, {\it Zap. Mineral. Imper. S. Petersburgskogo Ob$\check{s}$$\check{c}$},
{\bf 21}(2) (1885), 1-279.
\bibitem{furt}P. Furtw\"angler, \"Uber Gitter konstanter Dichte, {\it Monatsh. Math. Phys.} {\bf 43} (1936), 281-288.
\bibitem{grs}N. Gravin, S. Robins and D. Shiryaev, Translational tilings by a polytope, with multiplicity. {\it Combinatorica}
{\bf 32} (2012), 629-649.
\bibitem{gkrs}N. Gravin, M. N. Kolountzakis, S. Robins and D. Shiryaev, Structure results for multiple tilings in 3D.
{\it Discrete Comput. Geom.} {\bf 50} (2013), 1033-1050.
\bibitem{grub}P. M. Gruber and C. G. Lekkerkerker, {\it Geometry of Numbers} (2nd ed.), North-Holland, 1987.
\bibitem{hajo}G. Haj\'os, \"Uber einfache und mehrfache Bedeckung des $n$-dimensionalen Raumes mit einem W\"urfelgitter, {\it Math. Z.}
{\bf 47} (1941), 427-467.
\bibitem{hilb}D. Hilbert, Mathematical Problems, {\it G\"ottinger Nachr.} (1900), 253-297. {\it Bull. Amer. Math. Soc.} {\bf 8} (1902), 437-479;
{\bf 37} (2000), 407-436.
\bibitem{kolo}M. N. Kolountzakis, On the structure of multiple translational tilings by polygonal regions. {\it Discrete Comput. Geom}.
{\bf 23} (2000), 537-553.
\bibitem{lazo}J. C. Lagarias and C. Zong, Mysteries in packing regular tetrahedra, {\it Notices Amer. Math. Soc.}
{\bf 59} (2012), 1540-1549.
\bibitem{mann}C. Mann, J. McLoud-Mann and D. Von Derau, Convex pentagons that admit $i$-block transitive tilings, {\it Geom. Dedicata} (2018), in press.
\bibitem{mcmu}P. McMullen, Convex bodies which tiles space by translation, {\it Mathematika} {\bf 27} (1980), 113-121.
\bibitem{mink}H. Minkowski, Allgemeine Lehrs\"atze \"uber konvexen Polyeder, {\it Nachr. K. Ges. Wiss. G\"ottingen, Math.-Phys. KL}. (1897), 198-219.
\bibitem{rao}M. Rao, Exhaustive search of convex pentagons which tile the plane, arXiv:1708.00274
\bibitem{rein}K. Reinhardt, \"Uber die Zerlegung der Ebene in Polygone, {\it Dissertation,} Universit\"at Frankfurt am Main, 1918.
\bibitem{robi}R. M. Robinson, Multiple tilings of $n$-dimensional space by unit cubes, {\it Math. Z.} {\bf 166} (1979), 225-275.
\bibitem{stog}M. I. $\check{S}$togrin, Regular Dirichlet-Voronoi partitions for the second triclinic group (in Russian),
{\it Proc. Steklov. Inst. Math.} {\bf 123} (1975).
\bibitem{venk}B. A. Venkov, On a class of Euclidean polytopes, {\it Vestnik Leningrad Univ. Ser. Mat. Fiz. Him}. {\bf 9} (1954), 11-31.
\bibitem{voro}G. F. Voronoi,  Nouvelles applications des paramm\`etres continus \`a la th\'eorie des formes quadratiques.
Deuxi\`eme M\'emoire. Recherches sur les parall\'elo\`edres primitifs, {\it J. reine angew. Math.} {\bf 134} (1908), 198-287; {\bf 135} (1909), 67-181.
\bibitem{yz}Q. Yang and C. Zong, Multiple lattice tilings in Euclidean spaces, {\it Bull. London Math. Soc.} in press.
\bibitem{zong05}C. Zong, What is known about unit cubes. {\it Bull. Amer. Math. Soc.} {\bf 42} (2005), 181-211.
\bibitem{zong06}C. Zong, {\it The Cube: A Window to Convex and Discrete Geometry.} Cambridge University Press, Cambridge, 2006.
\bibitem{zong14}C. Zong, Packing, covering and tiling in two-dimensional spaces, {\it Expo. Math.} {\bf 32} (2014), 297-364.
\bibitem{zong17}C. Zong, Characterization of the two-dimensional five-fold lattice tiles, arXiv:1712.01122.
\end{thebibliography}

\vspace{0.6cm}
\noindent
Qi Yang, School of Mathematical Sciences, Peking University, Beijing 100871, China

\vspace{0.3cm}
\noindent
Corresponding author:

\noindent
Chuanming Zong, Center for Applied Mathematics, Tianjin University, Tianjin 300072, China

\noindent
Email: cmzong@math.pku.edu.cn

\end{document}